%
%
 \magnification \magstep1
\font\bold=cmbx10 at 14pt

\font\tts=cmtt8
\centerline{\bold Three-Dimensional Manifolds, Skew-Gorenstein Rings }
\medskip 
\centerline{\bold  and their Cohomology.}
\bigskip
\centerline{ Jan-Erik Roos}
\centerline{Department of Mathematics}
\centerline{Stockholm University}
\centerline{SE--106 91 Stockholm, SWEDEN}
\centerline{ e-mail: {\tt jeroos@math.su.se}}
\bigskip
 \rightline{ \it Dedicated to Ralf Fr{\"o}berg and Clas L{\"o}fwall at their $65^{th}$ birthdays.}
\bigskip

\centerline{May 12, 2010}

\bigskip

\def\mysec#1{\bigskip\centerline{\bf #1}\nobreak\par}

\def\cite#1{~[{\bf #1}]}
{\openup 1\jot

\mysec{ Abstract.}
\bigskip
Graded skew-commutative rings occur often in practice. Here are two examples:
1) The cohomology ring of a compact three-dimensional manifold.
2) The cohomology ring of the complement of a hyperplane arrangement (the Orlik-Solomon algebra).
We present some applications of the homological theory of these graded skew-commutative rings. In particular we find compact oriented 3-manifolds 
without boundary for which the Hilbert series of 
the Yoneda Ext-algebra of the cohomology ring of the fundamental group
is an explicit transcendental function.
This is only possible for large first Betti numbers of the 3-manifold (bigger than -- or maybe equal to -- 11).
We give also examples of 3-manifolds where the Ext-algebra of the cohomology ring of the fundamental group
is not finitely generated.


\mysec{0. Introduction}
\bigskip
Let $X$ be an oriented compact $3$-dimensional manifold without boundary.
The cohomology ring $H=H^*(X,{\bf Q})$
is a graded skew-commutative ring whose augmentation ideal $\bar H$ satisfies ${\bar H}^4=0$.
The triple (cup) product $x \cup y \cup z =\mu(x,y,z).e$,
where e is the orientation generator of $H^3$, defines a skew-symmetric trilinear form
on $H^1$ with values in $\bf Q$  (i.e. a {\it trivector}) and conversely, according to a theorem of Sullivan [Sul] 
any such form comes in this way from a 3-manifold $X$ (not unique)
whose cohomology algebra can be reconstructed from $\mu$ since by Poincar\'e duality $H^2 \simeq (H^1)^*$.
In the more precise case when $H^*$ is a also a Poincar\'e duality algebra, i.e.
the cup product $H^1 \times H^2 \rightarrow H^3$ is non-degenerate, it
follows that $H^*$ is a Gorenstein ring (cf. section 1 below).
 Such Gorenstein rings will be studied here. Any 3-manifold $M$ can be decomposed in a unique way
as a connected sum of prime 3-manifolds:
$$
M = P_1 \sharp P_2 \cdots \sharp P_k
$$
of prime manifolds $P_i$ (cf e.g. Milnor [Mil], Theorem 1), and with the exception of $S^3$ and $S^2 \times S^1$
any prime manifold is also irreducible ([Mil] and for them $\pi_2(M)=0$ ([Mil], Theorem 2). If furthermore
$\pi_1(M)$ is infinite then also the higher homotopy groups are 0 ($M$ is said to be aspherical)
so that $M$ is the Eilenberg-Maclane space
$K(\pi_1(M),1)$ and the cohomology ring of $M$ is isomorphic to the cohomology ring of the group 
$\pi_1(M)$  (this is of course also true for any 3-manifold $X$ which is aspherical).
(For most of the applications below we suppose that $H=H^*(X)$ is a Poincar\'e
duality algebra and we suppose that the base field $k$ is of characteristic 0, the preference being {\bf Q}.)
The ring $H$ has interesting homological properties
which have not yet been fully studied, and we wish to continue such a study here.
For small values of the first Betti number $b_1(X)= {\rm dim}_{\bf Q}H^1(X,{\bf Q})$
 the ring $H$ is a Koszul algebra (cf. section 2 below) so that in particular the generating series
$$
P_H(z) = \sum_{i\geq 0}|{\rm Tor}_i^H(k,k)|z^i = H(-z)^{-1}  \leqno(0.1)
$$
where $H(z) = 1 + |H^1(X,Q)|z + |H^2(X,Q)|z^2+z^3$,
and where $|V|$ denotes the dimension of the vector space $V$. 
 But for bigger Betti numbers many new
phenomena occur. In particular we will see that for $b_1(X)= 12$
(and maybe even for $b_1(X) =11$) there are a few  examples where
$P_H(z)$ is an explicit transcendental function (thus we are far away
from the Koszul case of formula (0.1)!). However, maybe 11 is the best possible number here.
On the other hand for bigger $b_1(X)$ the possible series $P_H(z)$ are
rationally related to the family of series which occurred in
connection with the Kaplansky-Serre questions a long time ago [An-Gu].
But even for smaller $b_1(X)$ other strange homological properties of $H$  
occur: we will give examples (probably best possible) where $b_1(X)=11$ and the Yoneda ${\rm Ext}$-algebra
${\rm Ext}^*_H(k,k)$ is not finitely generated.  
The implications of all this for $3$-manifold groups have not been
fully explored. Cf. also [Sik]. Conversely, the connection with $3$-manifolds
makes it possible to go backwards and to deduce results in the homology theory of skew-commutative algebras (and there are also related
results in the {\it commutative} case).
 Note that we are studying everything in characteristic $0$.
There seem to be some relations with Benson [Be] but he works
over finite fields.
Finally, let us remark that even for the special 3-dimensional case when $X$
is the boundary manifold of a line arrangement in $P^2({\bf C})$ 
we can have e.g. the same strange transcendental phenomenon as above but the prize to pay for this
is to accept even bigger $b_1(X)$. The three-manifolds that occur here are called
graphic manifolds [Co-Su] and they are aspherical if the line arrangement is not a pencil 
of lines. 
\bigskip
\mysec{1. Graded skew-Gorenstein rings and classification of trivectors}
\bigskip
Let us first recall that a local commutative Gorenstein ring was defined in [Bass] as
a ring $R$ which has a finite injective resolution as a module over itself.
In particular if $R$ is artinian this means that $R$ is injective as a module over itself.
It is also equivalent to saying that the socle of $R$ , i.e. ${\rm Hom}_R(R/m,R)$ ($m$ is the maximal ideal of $R$)
is 1-dimensional over $k=R/m$. 
Things are more complicated in the noncommutative case [Fo-Gr-Rei], but if
$R$ is skew-commutative artinian the Gorenstein property is equivalent to $R$ being injective as a module
over itself (left or right --- these two conditions are equivalent --- and they are also equivalent to
saying that the left --- or right ---  socle of $R$ is one-dimensional).   
In the special case when $X$ is an oriented compact three-dimensional manifold without boundary, and
when  $H = H^*(X,{\bf Q})$ is the cohomology ring of $X$, we let $R=H=H^*(X,{\bf Q})$.
 From the preceding definition it follows that
$R$ is Gorenstein if and only if $R$ is generated by $H^1$ and $H$ is a Poincar\'e duality algebra
(we assume that $|H^1| > 1$). 

We now turn to the classification of such $R$:s when $|H^1(X,{\bf Q})|\leq 8$. We will use the
classification of trivectors in $H^1$ described in section 35 of the book [Gur].
For a background we refer the reader to our introduction, to [Sul] and to the section 4 [Intermezzo: Classification of skew-symmetric forms]
of [Sik]. Since we are only studying non-degenerate trivectors, their ranks are equal to
the dimension of $H^1$.

For the rank 5 there is only one trivector (denoted by III in [Gur], page 391) and which can be denoted
by (here $e_1,e_2,e_3,e_4,e_5$ is a basis for $H^1$ and the $e^i$:s are the dual basis elements in $(H^1)^*$):
$$
e^1\wedge e^4 \wedge e^3 + e^2 \wedge e^5 \wedge e^3
$$
and the corresponding Gorenstein ring
is a Koszul algebra (the preceding trivector is the correct one as in [Gur] -- in [Sik] there
is a minor misprint in its form). 
It is only when the ranks are $\geq 6$ that non-Koszul Gorenstein rings occur.
Let us give the details in the first nontrivial case of [Gur], namely case IV of rank 6, page 391:
$$
f=e^1 \wedge e^2 \wedge e^3 +e^3 \wedge e^4 \wedge e^5 + e^2 \wedge e^5 \wedge e^6  
$$
Note that $H^*$ is a quotient of the exterior algebra 
$$
E(e_1,e_2,e_3,e_4,e_5,e_6)
$$
by an ideal that we want to determine. To begin with we want to determine 
 all quadratic elements $g=\sum_{j<k}g_{j,k}e_j\wedge e_k$ that go to 0 in the quotient $H^2$. Since the the pairing $H^1 \otimes H^2 \rightarrow H^3$ is nondegenerate this boils down to
determine when $e_s\wedge g$ goes to zero for $s=1,\ldots,6$. Using $f$ we get the conditions
$\sum_{j<k}g_{j,k}f(e_s,e_j,e_k) = 0$
for $s=1,\ldots, 6$.
This gives, using the explicit form of $f$ and the condition that $f$ is skew-symmetric (note that if $f$ were only the monomial $f_{mon}=e^1\wedge e^2 \wedge e^3$, then
$f_{mon}(e_{i_1},e_{i_2},e_{i_3})$ is non zero (=$\pm 1$) if and only if $i_1,i_2,i_3$ is a permutation of $1,2,3$)
 the six conditions
$s = 1,\ldots,6$:
$$ 
g_{2,3} =0,\quad -g_{1,3}+g_{5,6}=0,\quad g_{1,2}+g_{4,5}=0,\quad g_{3,5} =0,\quad g_{3,4}-g_{2,6}=0,\quad
g_{2,5}=0.
$$
The 9 solutions of this system of 6 linear equations
lead to the relations $e_1\wedge e_3+e_5\wedge e_6$,$e_1\wedge e_2-e_4\wedge e_5$ and
 $e_3\wedge e_4+e_2\wedge e_6$ and to 6 monomial relations, 
leading to the following ring with quadratic relations (we do not write wedge for
multiplication):
$$
R_{IV} = {E(e_1,e_2,e_3,e_4,e_5,e_6) \over (e_1e_2-e_4e_5,e_1e_3+e_5e_6,e_2e_6+e_3e_4,e_1e_4,e_1e_5,e_1e_6,e_2e_4,e_3e_6,e_4e_6)} 
$$
But this quotient ring $R_{IV}$ has Hilbert series $R_{IV}(z)=1+6z+6z^2+2z^3$. This means that we have to find a last
cubic relation. One finds that the cube of the maximal ideal of $R_{IV}$
is generated by $e_2e_5e_6$ and $e_2e_3e_5$.
But $f(e_2,e_3,e_5)=0$ and $f(e_2,e_5,e_6)=1$ so that the corresponding 
Gorenstein ring is $G_{IV}=R_{IV}/e_2e_3e_5$ and we will see that this Gorenstein ring is {\it not} a Koszul algebra.
For the case V of {\it loc.cit.} the $f$ is given by
$$
f = e^1 \wedge e^2 \wedge e^3 + e^4 \wedge e^5 \wedge e^6
$$
leading in the same way to the ring with quadratic relations:
$$
R_V ={E(e_1,e_2,e_3,e_4,e_5,e_6) \over(e_1e_4,e_1e_5,e_1e_6,e_2e_4,e_2e_5,e_2e_6,e_3e_4,e_3e_5,e_3e_6)}
$$
which also has Hilbert series $R_{V}(z)=1+6z+6z^2+2z^3$, but in this case we have to divide by
$e_1e_2e_3-e_4e_5e_6$ to get the Gorenstein ring $G_V=R_V/(e_1e_2e_3-e_4e_5e_6)$
which is, as we will see below, not a Koszul algebra either.
\bigskip
In [Gur, pages 393-395] the classification of 3-forms of rank 7 is given as the 5 cases VI, VII, VIII, IX, X and those forms
of rank 8  are described as the 13 cases XI, XII,..., XXIII.
The classification of 3-forms of rank 9 are given in [Vin-El].
We will describe the homological behaviour of the corresponding Gorenstein rings in the section 4.
 \mysec{2. Calculating the Koszul dual of the Gorenstein ring associated to a 3-form.}
\bigskip
Let $R$ be any finitely presented ring (connected $k$-algebra) generated in degree 1 and having quadratic relations.
It can be described as the quotient $T(V)/(W)$, where $T(V)$ is the tensor algebra on a (finite-dimensional)
$k$-vector space $V$, placed in degree 1 and $(W)$ is the ideal in $T(V)$, defined 
by a sub-vector space $W \subset V\otimes_k V$. The Yoneda Ext-algebra of R is defined by
$$
{\rm Ext}^*_R(k,k) = \bigoplus_{i\geq 0}{\rm Ext}^i_R(k,k) 
$$ 
where $k$ is an $R$-module in the natural way and where the multiplication is the Yoneda product.
The sub-algebra of ${\rm Ext}^*_R(k,k)$ generated by ${\rm Ext}^1_R(k,k)$ is called the Koszul dual of $R$ and it is
denoted by $R^!$. It can be calculated as follows: consider the inclusion map $W \rightarrow  V\otimes_k V$.
Taking k-vector space duals (denoted by $W^*$ and $(V\otimes_k V)^*$ we get the exact sequence
($W^{\perp}$ = those linear $f: V\otimes_k V \rightarrow k$ that are $0$ on $W$):
$$
   0 \leftarrow W^* \leftarrow (V\otimes_k V)^* \leftarrow W^{\perp} \leftarrow 0
$$ 
Now $R^! = T(V^*)/(W^{\perp})$  (we have used that  $(V\otimes_k V)^* =V^*\otimes_k V^*$).

For all this cf. [L\"of 2]. Note that $(R^!)^!$ is isomorphic to $R$.
Note also that if $R$ also has cubic relations and/or higher relations then
${\rm Ext}^*_R(k,k)$ is still defined, and the subalgebra generated by ${\rm Ext}^1_R(k,k)$ is still given by the formula
 $T(V^*)/(W^{\perp})$ where $W$ is now only the ``quadratic part'' of the relations of $R$. In particular
$(R^!)^!$ is now only isomorphic to $T(V)$ divided by the quadratic part of the relations.
Note also that in general, if $R$ is skewcommutative then ${\rm Ext}^*_R(k,k)$ is a cocommutative Hopf algebra
which is the enveloping algebra of a graded Lie algebra, and $R^!$ is a sub Hopf algebra which is the
enveloping algebra of a smaller graded Lie algebra. (All this is also true if $R$ is commutative,
but now the Lie algebras are super Lie algebras.)
In [L\"of 2] Corollary 1.3, pages 301-302 there is a recipe about how to calculate the Koszul dual of an
algebra with quadratic relations.
Applying this to the case of $R_{IV}$ of the previous section we find that $R_{IV}^!$ is the algebra
$$
{k<X_1,X_2,X_3,X_4,X_5,X_6> \over
([X_1,X_2]+[X_4,X_5],[X_3,X_4]-[X_2,X_6],[X_1,X_3]-[X_5,X_6],[X_2,X_3],[X_2,X_5],[X_3,X_5])}
$$
where $[X_i,X_j]=X_iX_j-X_jX_i$ is for $i<j$ the Lie commutator of $X_i$ and
$X_j$.
In general if one starts with a (skew)-commutative algebra $A$ it is often easy to calculate the Hilbert series
of $A$. But calculating the Hilbert series of $A^!$ is often very difficult and this series can even be a
transcendental function.
But in this case it is rather easy: we get that $R_{IV}^!(z) = 1/(1-6z+6z^2-2z^3)$ so that 
$R_{IV}(-z)R_{IV}^!(z) = 1$ and we even have (the in general strictly stronger assertion ([Ro 4], [Pos]))
that $R_{IV}$ is a Koszul algebra.
(For the definitions and equivalent characterizations of Koszul algebras we refer to
[L\"of 2, p. 305, Theorem 1.2].)  A similar result for $R_V$ holds true, and in this case we
can directly apply a result of Fr\"oberg [Fr\"o 1], since $R_V$ has quadratic monomial relations. 
But neither the Gorenstein quotients $G_{IV}$ nor $G_V$ are Koszul algebras and in the next section
we will see how to relate their Hilbert series and the corresponding Hilbert series of their Koszul duals
to their two-variable Poincar\'e-Betti series 
$$
P_G(x,y)=\sum_{i,j} |{\rm Tor}_{i,j}^G(k,k)|x^iy^j
$$
But we will first deduce a few results about explicitly calculating the Gorenstein ring and its 
Koszul dual associated to a given 3-form.

An old result of Macaulay gives a nice correspondence between commutative artinian graded Gorenstein algebras having
socle of degree $j$
of the form $k[x_1,x_2,\ldots,x_n]/I$ and homogeneous forms of degree $j$ in the dual of $k[x_1,x_2,\ldots,x_n]$
({\it cf. eg.} Lemma 2.4 of [Elk-Sri]). Here is a skew-commutative version 
(here described only for socle  degree 3), the proof of which follows from
Exercise 21.1, page 547 of [Eis], where we have to replace polynomial rings by exterior algebras 
(I thank Antony Iarrobino who suggested that such a result should be true):

PROPOSITION 2.1 (``Skew''-Macaulay).- Let $T=E[y_1,\ldots, y_n]$ be the exterior algebra in $n$ variables of degree 1
over a field of characteristic $0$. Consider a polynomial skew-differential operator $D$ on $T$ with constant
coefficients of the form  
$$
D =\sum_{i_1,\ldots,i_n} a_{i_1,\ldots,i_n} (\partial /\partial y_1)^{i_1}\cdots (\partial /\partial y_n)^{i_n}
$$
(the $i_j$ are 0 or 1). The symbol of $D$ is the $\sum_{i_1,\ldots,i_n} a_{i_1,\ldots,i_n}x_1^{i_1}\cdots x_n^{i_n}
$ which is obtained by replacing each $\partial /\partial y_i$ by $x_i$ in a new exterior algebra
$S=E[x_1,\ldots,x_n]$. Let now $0 \neq f \in T$ be a homogeneous polynomial of degree 3. Let $I$ be the set of 
symbols in $S$ of differential operators $D$ as above, such that $Df=0$. Then $S/I$ is a zero-dimensional
Gorenstein ring and there is a converse assertion.

We will not use the preceding result, which is only included here for historical reasons.
\bigskip
We now go back to the reasoning that we used when we treated the case IV in section 1.
 Recall first that in ([L\"of 2], Corollary 1.3, p. 301) it is proved in particular that if we
have a quadratic algebra $R=E[x_1,\ldots,x_n]/(f_1,\ldots,f_r)$ where the 
$$
f_i=\sum_{j < k}b_{i,j,k}x_j\wedge x_k, \,\,b_{i,j,k} \in k \quad, i=1,\dots r,
$$ 
then 
$$
R^! = k<X_1,\ldots,X_n> /({\phi}_1,\ldots,{\phi}_s),
$$
where
$$
{\phi}_i = \sum_{j < k}c_{i,j,k}[X_j,X_k], \,\, c_{i,j,k} \in k,\quad i = 1,\ldots,s,
$$
where $[X_j,X_k]=X_jX_k-X_kX_j$ and $(c_{i,j,k})_{jk},\quad i=1,\ldots s$, is a basis of the solutions of system of linear equations:
$$
\sum_{j<k}b_{i,j,k}X_{jk},\quad i = 1,\ldots,r\quad ({\rm hence}),\quad s = {n+1 \choose 2}-r.
$$
Now given a skew-symmetric 3-form $\Psi(e^1,e^2,\ldots,e^n)$ of rank $n$, written as a linear combination of 
$e^{i_1}\wedge e^{i_2}\wedge e^{i_3}$, where $i_1 <i_2<i_3$ (with coefficients that are $a_{i_1i_2i_3}$),
 we have that 
$$
f_i = \sum_{j<k}b_{ijk} e_j\wedge e_k, \quad i=1,\ldots,r
$$
 are relations corresponding to $\Psi$ (cf. the discussion of case IV in section 1) if and only if
$$
\sum_{j<k}b_{i,j,k}\Psi(e_i,e_j,e_k) = 0,\quad {\rm for}\,\, i=1,\ldots, n.
$$
But for fixed $i$, $\Psi(e_i,e_j,e_k)$ is $a_{ijk}$ if and only if $i<j<k$
and $j<k<i$, and  $-a_{ijk}$ if and only if $j<i<k$ . This corresponds to taking the graded skew-derivative of $\Psi$ with
respect to $e^i$. 
Combining this with the L\"ofwall description of the calculation of $R^!$ given above,
we finally arrive at the following useful Theorem-Recipe to calculate the Koszul dual $G^!$ of the Gorenstein
ring associated to a skew 3-form:
\bigskip   
THEOREM-RECIPE 2.1 --- Let $\Psi(e^1,e^2,\ldots,e^n)$ be a skew-symmetric 3-form of rank $n$, $X$ one of the 3-manifolds
giving rize to $\Psi$ and $G = H^*(X,{\bf Q})$. The Koszul dual $G^!$ of $G$ is obtained as follows:
Calculate the $n$ skew-derivaties of $\Psi$ with respect to the variables $e^1,e^2,\ldots,e^n$.
Then
$$
G^! \simeq {k<X_1,X_2,\ldots,X_n>\over (q_1,q_2,\ldots,q_n)}
$$
where $k<X_1,X_2,\ldots,X_n>$ is the free associative algebra generated by the variables $X_i$
that are {\it dual} to the $e^i$:s in $\Psi$ and where the $q_i$:s are obtained by replacing
each quadratic element $e^s\wedge e^t$ ($s < t$) in the $\partial \Psi /\partial e^i$ by the commutator
$[X_s,X_t] =X_sX_t-X_tX_s$ in $k<X_1,X_2,\ldots,X_n>$, for $i = 1,\ldots n$.

Note that if $G^!$ is given it is easy to calculate ``backwards'' the ring with quadratic relations
$(G^!)^!$. After that it is easy to find the extra cubic relations we should divide with to get $G$.
Using this theorem-recipe we finally find:
\bigskip
THEOREM 2.2 --- a) The double Koszul duals $(G^!)^!$ of the Gorenstein rings corresponding trivectors
of rank $= 7$ are Koszul algebras in the cases VI, VII, VIII, IX  and X and therefore $G^!$ is a Koszul
algebra in cases VI, VII, VIII, IX, X. (But they are not isomorphic). The corresponding Gorenstein algebras
$G_{VI},G_{VII},G_{VIII},G_{IX},G_{X}$ are also Koszul algebras.
\bigskip
b) When it comes to trivectors of rank 8, i.e. the cases
XI, XII, $\ldots$, XXIII the situation is more complicated already from the
homological point of view:
Indeed in case XI the double Koszul dual $(G^!)^!$ is already a 
Gorenstein ring with quadratic relations, thus equal to $G$ and $G^!$ is
{\it not} a Koszul algebra.
But the case XII treated explicitly below is slightly different and not a Koszul algebra either. 
The cases XIII, XIV and XV are Koszul algebras.
But case XVI is as above: $(G^!)^!$ is a Koszul algebra, but we have
to  divide out a cubic form to get $G_{XVI}$ which is therefore not Koszul.
Finally the algebras corresponding to the cases XVII, XVIII, XIX, XX, XXI, XXII and XXIII are Koszul algebras (but not isomorphic). 
\bigskip 
EXAMPLE 2.1 Here is a use of the THEOREM-RECIPE 2.1:
Let us consider the case XII of rank 8.
Here the 3-form $f_{XII}=\Psi$ is
$$
\Psi(e^1,e^2,\ldots,e^8)=e^5\wedge e^6 \wedge e^7 +e^1 \wedge e^5 \wedge e^4 + e^2\wedge e^6 \wedge e^4 
+e^3\wedge e^7 \wedge e^4+ e^3\wedge e^6 \wedge e^8 
$$ 
We calculate the 8 partial derivaties of this skew form (we do not write out the $\wedge$ sign):
$$
{\scriptstyle \partial \Psi /\partial e^1 = e^5e^4,\quad \partial \Psi /\partial e^2 = e^6e^4, 
\quad \partial \Psi /\partial e^3 = e^7e^4+e^6e^8,\quad \partial \Psi /\partial e^4 = e^1e^5+e^2e^6+e^3e^7,}
$$
$$
{\scriptstyle \partial \Psi /\partial e^5 = e^6e^7-e^1e^4,\quad \partial \Psi /\partial e^6 = -e^5e^7-e^2e^4-e^3e^8,\quad
\partial \Psi /\partial e^7 = e^5e^6-e^3e^4, \quad \partial \Psi /\partial e^8 = e^3e^6} 
$$
leading to $G^! = k<X_1,X_2,X_3,X_4,X_5,X_6,X_7,X_8>$ divided by the ideal
$$
 \scriptstyle {( [X_5,X_4],[X_6,X_4],[X_7,X_4]+[X_6,X_8],
[X_1,X_5]+[X_2,X_6]+[X_3,X_7],[X_6,X_7]-[X_1,X_4],}
$$
$$
 \scriptstyle {
[X_5,X_7]+[X_2,X_4]+[X_3,X_8],[X_5,X_6]-[X_3,X_4],[X_3,X_6])}
$$
and the Hilbert series $1/G^!(z)= 1-8z+8z^2-z^3-z^4$. Furthermore one sees that $(G^!)^!$ has Hilbert series
 $1+8z+8z^2+z^3$  so that we do not have to divide by cubic elements to get $G_{XII}$ which is non-Koszul.
As a matter of fact we have a generalization of a formula of L\"ofwall (cf. next section), giving that
$$
{1\over P_{G_{XII}}(x,y)} = (1+1/x)/G_{XII}^!(xy)-G_{XII}(-xy)/x
$$
so that 
$$
{1\over P_{G_{XII}}(x,y)} = 1-8xy+8x^2y^2-x^3y^3-x^3y^4-x^4y^4
$$
leading to {\it e.g.}  ${\rm Tor}_{3,4}^{G_{XII}}(k,k)$ being 1-dimensional.
The other cases I-XXIII are treated in a similar manner: sometimes $G^!(z)$
can be directly calculated since we have a finite Groebner basis for the
non-commutative ideal and in all cases using the Backelin et al programme BERGMAN [Ba].
\mysec{3. A generalization of a formula of L\"ofwall to Gorenstein rings with $m^4=0$}
\bigskip
In his thesis [L\"of 2] Clas L\"ofwall proved in particular that if $A$ is any graded connected
algebra with ${\bar A}^3 =0$ then the double Poincar\'e-Betti series is given by the formula
$$
{1\over P_A(x,y)} = (1+1/x)/A^!(xy)-A(-xy)/x , \leqno(3.1)
$$
In [Ro 1] we have an easy proof of this in the case when $A$ is commutative (works also 
in the skew-commutative case [Ro 2]) using the fact that in these cases $A^!$ is a sub-Hopf algebra of the 
big Hopf algebra ${\rm Ext}^*_A(k,k)$ and according to a theorem of Milnor-Moore this big Hopf algebra is
free as a module over $A^!$. We now show
\bigskip
THEOREM 3.1.--- The formula (3.1) is true when $(R,m)$ is a graded (skew-) commutative Gorenstein ring
with $m^4=0$.
\medskip
PROOF: First we note that Avramov-Levin have proved [Av-Le] in the commutative case that the natural map
$$
 R \rightarrow R/soc(R)
$$
is a Golod map (cf. section 5 for the skew-commutative case that we use here).
 Now the socle of $R$ is $m^3$ which is $s^{-3}k$
Therefore we have
$$
P_R/soc(R)(x,y) = {P_R(x,y)\over 1-x(P_R^{R/soc(R)}(x,y)-1)} = {P_R(x,y)\over 1-x^2y^3P_R(x,y)} \leqno(3.2)
$$
But $R/m^3$ is a local ring where the cube of the maximal ideal is 0.
Thus the formula of L\"ofwall can be applied, and we obtain using that $(R/m^3)^! = R^!$ the formula
$$
{1\over P_R/m^3(x,y)} = (1+1/x)/R^!(xy)-(R/m^3)(-xy)/x
$$
which combined with (3.2) gives (3.1) since $R(z) = R/m^3(z) +z^3$.
\bigskip
REMARK 3.1.--- Thus it is clear that we can not obtain cases where the ring
$H^*(X,{\bf Q})$ has ``bad'' homological properties if $b_1(X) \leq 8$.
We therefore study the case when $b_1(X)=9$. In this case there is a classification of the tri-vectors
by E.B. Vinberg and A.G. Elashvili [Vin-El].
But even in this case it seems impossible to get ``exotic'' $H^*(X,{\bf Q})$.
In fact we have e.g. by the procedure above analyzed all cases
in the Table 6 (pages 69-72) of [Vin-El] (those cases that have a * in loc.cit. correspond to $b_1(X)<9$ and
they have already been treated). 
\bigskip
THEOREM 3.2.--- For the 3-forms of rank 9 of Table 6 (pages 69-72) of [Vin-El] we have that the corresponding
$1/R^!(z)=1-9z+9z^2-z^3$, and the corresponding Gorenstein ring is a Koszul algebra in all cases except
\medskip
a) the cases 79, 81, 85 where $1/R^!(z)= 1-9z+9z^2-3z^3$
\medskip
b) the case 83    where $1/R^!(z)= 1-9z+9z^2-z^3-5z^4+4z^5-z^6$
\bigskip
PROOF.- The proof is by using the THEOREM-RECIPE 2.1  above. For the Koszul cases we refer the reader to the
Appendix. Here
we only indicate what happens in the exceptional cases a) and b).
In case 79 in a) the 3-form is (we do not write out the $wedge$ for multiplication):
$$
f_{79}=e^1e^2e^9+e^1e^3e^8+e^2e^3e^7+e^4e^5e^6
$$
We take the 9 partial derivatives
$\partial f_{79} /\partial e^i$ for $i=1\ldots 9$
and we obtain:
$$
e^2e^9+e^3e^8,\quad -e^1e^9+e^3e^7,\quad -e^1e^8-e^2e^7,\quad e^5e^6,\quad -e^4e^6,\quad e^4e^5,
\quad e^2e^3,\quad e^1e^3,\quad e^1e^2
$$
leading to $R^! = $
$$
 {{\scriptstyle k<X_1,X_2,X_3,X_4,X_5,X_6,X_7,X_8,X_9>} \over {\scriptscriptstyle ( [X_2,X_9]+[X_3,X_8],
-[X_1,X_9]+[X_3,X_7],[X_1,X_8]+[X_2,X_7],[X_5,X_6],[X_4,X_6],[X_4,X_5],[X_2,X_3],[X_1,X_3],[X_1,X_2])}}
$$
Now it turns out that the Gr\"obner basis of the ideal above is finite
and in degree $\leq 2$. This proves that $R^!$ is Koszul and has
$R^!(z) = 1-9z+9z^3-3z^3$, since the commutative ring $(R^!)^!$ has
Hilbert series $1+9z+9z^2+3z^3$.
The ideal $m^3$ is generated by the three elements ($e_1e_2e_9,e_4e_5e_6,e_1e_2e_3$)
and
$$
f_{79}(e_1,e_2,e_9)=1,\quad f_{79}(e_4,e_5,e_6)=1 \quad {\rm and}\quad f_{79}(e_1,e_2,e_3)=0
$$ 
It follows that if we divide out by the two extra elements $e_1e_2e_9-e_4e_5e_6,e_1e_2e_3$ we do have the non-Koszul
Gorenstein ring $G_{79}$ we are looking for. The rings $G_{81}$ and $G_{85}$ in a) are treated in the same way.
They correspond to:
$$
f_{81}= e^1e^2e^9+e^1e^3e^8+e^1e^4e^6+e^2e^3e^7+e^2e^4e^5+e^3e^5e^6
$$
and
$$
f_{85}=e^1e^2e^9+e^1e^3e^5+e^1e^4e^6+e^2e^3e^7+e^2e^4e^8
$$ 
\medskip
In case 83 of b) the form 
$$
f_{83}=e^1e^2e^9+e^1e^3e^5+e^1e^4e^6+e^2e^3e^7+e^2e^4e^8+e^3e^4e^9 
$$
and its 9  skew-derivaties lead to the quotient:
$$
\scriptstyle {k<X_1,X_2,X_3,X_4,X_5,X_6,X_7,X_8,X_9>}
$$
divided by the ideal
$$
\scriptstyle { ( [X_2,X_9]+[X_3,X_5]+[X_4,X_6],
-[X_1,X_9]+[X_3,X_7]+[X_4,X_8],-[X_1,X_5]-[X_2,X_7]+[X_4,X_9],}
$$
$$\scriptstyle {[X_1,X_6]+[X_2,X_8]+[X_3,X_9],
[X_1,X_3],[X_1,X_4],[X_2,X_3],[X_2,X_4],[X_1,X_2]+[X_3,X_4])}
$$
In this case the Hilbert series of $R^!$ is
$1/(1-9z+9z^2-z^3-5z^4+4z^5-z^6)$,
and the ring $(R^!)^!$ has Hilbert series:
$1+9z+9z^2+z^3$ so that $G_{83} = (R^!)^!$ which is a Gorenstein non-Koszul algebra and the THEOREM 3.2 is proved.
\bigskip
 
Thus to obtain examples of $X$ where the homological properties of the cohomology algebra $H^*(X,Q)$ are complicated
we need to study the cases where the dimension of $H^1(X,{\bf Q})$ is $\geq 10$.
But the trivectors of rank 10 have not been classified, and even the cases
of Sikora [Sik], Theorem 1 which uses for any simple Lie algebra $g$ 
the trilinear skew-symmetric form $\Psi_g(x,y,z) =\kappa(x,[y,z])$ where $\kappa$ is the Killing form of $g$,
do not seem to give anything exotic from our point of view, at least when the dimension of the Lie algebra is 10.
 We therefore need a new way of constructing 
Gorenstein rings corresponding to trivectors of rank $\geq 10$, but defined
in another way.
This will be done in the next section.
\mysec{4. New skew-Gorenstein rings.}
Let $(R,m)$ be any ring with $m^3 = 0$ which is the quotient of the
exterior algebra $E(x_1,\ldots,x_n)$ by homogeneous forms of degrees $\geq 2$
and let $I(k)$ be  the injective envelope of the residue field $k=R/m$ of $R$.
Then $G=R \propto I(k)$ ([Fo-Gr-Rei] is a skew-Gorenstein ring with maximal ideal
$n = m \oplus I(k)$ satisfying $n^4=0$, which therefore corresponds to a 
trivector of rank $|m/m^2|+|I(k)/mI(k)|=|m/m^2|+|m^2|$ and
according to the Sullivan theory it comes from a 3-manifold $X$ with
the dimension of $H^1(X,{\bf Q})$ being equal to $|m/m^2|+|m^2|$,
 It turns out that the homological properties of $G$ are closely
related to those of the smaller ring $R$.
Indeed we have the following general result:
\bigskip
THEOREM 4.1 (Gulliksen [Gu]).--- Let $R$ be (skew-)commutative and $R \propto M$ be the trivial extension of $R$ with the $R$-module $M$.
Then we have an exact sequence of Hopf algebras:
$$
k \rightarrow T(s^{-1}{\rm Ext}^*_R(M,k)) \rightarrow {\rm Ext}^*_{R \propto M}(k,k) \rightarrow {\rm Ext}^*_R(k,k) \rightarrow k
$$
where $T(s^{-1}{\rm Ext}^*_R(M,k))$ is the free algebra on the graded vector space $s^{-1}{\rm Ext}^*_R(M,k)$ and where
the arrow ${\rm Ext}^*_{R \propto M}(k,k) \rightarrow {\rm Ext}^*_R(k,k)$ is a split epi-morphism (there is a splitting
ring map $R\propto M \rightarrow R$).
In particular we have the formula for Poincar\'e-Betti series
$$
P_{R\propto M}(x,y) = {P_R(x,y)\over 1-xP_R^M(x,y)}
$$
\bigskip
REMARK 4.0.- Gulliksen's proof [Gu] is in the commutative setting and uses Massey operations.
In [L\"of 1], Clas L\"ofwall studies the more general case of a trivial extension $R\propto M$, where $R$
is any ring (not necessarily commutative) and $M$ is an $R$-bimodule.
 In section 6, pages 305-306 of [L\"of 1]  he derives the
Gulliksen formula when $R$ is commutative using an idea of mine (cf. [L\"of 1], page 288). In the same
way one deduces the Gulliksen formula in the skew-commutative case.
\bigskip
Now turn to the artinian case and
 assume that $M$ is finitely generated and  let $\tilde M = {\rm Hom}_R(M,I(k))$
be the Matlis dual of $M$ [Mat] (if $M$ as an $R$-module is a vector space over $k$, then this is   
the ordinary vector space dual)
In this case we have a well-known formula (cf e.g. Lescot [Les 2], Lemme 1.1) 
$$
 {\rm Ext}^*_R(M,k) \simeq  {\rm Ext}^*_R(k,\tilde M)
$$
In particular if $M = I(k)$ then $\tilde M ={\rm Hom}_R(I(k),I(k)) = R$ [Mat]
so that 
$$
 {\rm Ext}^*_R(I(k),k)\simeq  {\rm Ext}^*_R(k,R)
$$
Thus, if we use the Theorem on $R \propto I(k)$ we are led to the study of the Bass series of $R$,
i.e. the generating series in one or two variables of ${\rm Ext}^*_R(k,R)$.
It turns out that in many (most) of the cases we study here, the Bass series of $R$ divided by
the Poincar\'e-Betti series $P_R$ is a very nice explicit polynomial (for more details
about this -- the B{\o}gvad formula - we refer to section  5 of this paper).
Thus if we want strange homological properties of the Ext-algebra of the Gorenstein ring
$R \propto I(k)$ i.e the Ext-algebra of the cohomology ring of the corresponding fundamental group
$\pi_1(X)$ of the 3-manifold $X$ corresponding to the Gorenstein ring $R \propto I(k)$
we only have to find $(R,m)$ with $m^3=0$ with strange properties.
\bigskip
COROLLARY 4.1 ---  
In the case when $R$ is the cohomology ring (over ${\bf Q}$) of the complement
of a line arrangement $L$ in $P^2({\bf C})$ i.e. the Orlik-Solomon
algebra of $L$, the 3-manifold $X$ correponding to the
Gorenstein ring $R\propto I(R/m)$ can be chosen as the boundary manifold
of a tubular neighborhood of $L$ in $P^2({\bf C})$ [Co-Su].
In this case we found in [Ro 2] an arrangement where the Ext-algebra of $R$
was not finitely generated. Then the Ext-algebra of the cohomology ring of X that corresponds to
$R \propto I(k)$ can not be finitely generated either since by THEOREM 4.1  
it is mapped onto ${\rm Ext}^*_R(k,k)$.
In this case the dimension of $H^1(X,{\bf Q})$ is 12.

In [Ro 2] we found two cases of arrangements: the MacLane arrangement
and the mleas arrangement where the corresponding $R^!(z)$ is a transcendental function, and this leads
to two cases where $H^1(X,{\bf Q})$ is of dimension 20, resp. 21. 
In order to press down this dimension to 12 and maybe to 11, I have to use
some of my earlier results.
In our paper [Ro 3] describing the homological properties of quotients
of exterior algebras in 5 variables by quadratic forms (there are 49 cases found), we have
found 3 cases (cases 12, 15 and 20) where $R^!(z)$ is proved to be transcendental and is
explicitly given, and 3 other cases (cases 21, 22 and 33) where we conjecture that $R^!(z)$
is transcendental, but no explicit formula can be given, even in the case
33 where we have now calculated the series $R^!(z)$ up to degree 33 using Backelin's et al
programme BERGMAN (some details are given in [Ro 2], where the ``educated guess'' now has
to be abandoned):

Here is the Case 20:
$$
R_{20}={E(x_1,x_2,x_3,x_4,x_5)\over(x_1x_4+x_2x_3,x_1x_5+x_2x_4,x_2x_5+x_3x_4)}
$$
Here the Hilbert series is $R_{20}(z)=1+5z+7z^2$,
and the Koszul dual $R_{20}^!$ is given (according to the recipe we have described above) by
$$
k<X_1,X_2,X_3,X_4,X_5> \over {\scriptstyle([X_1,X_2],[X_1,X_3],[X_3,X_5],[X_4,X_5],[X_1,X_4]-[X_2,X_3],[X_1,X_5]-[X_2,X_4],[X_2,X_5]-[X_3,\
X_4])} \leqno(4.1)
$$
where $ k<X_1,X_2,X_3,X_4,X_5>$ is the free associative algebra in the five variables $X_i$ and
$[X_i,X_j]=X_iX_j-X_jX_i$ is the commutator.
The corresponding Hilbert series is:
$$
{1\over R_{20}^!(z)}=\prod_{n=1}^\infty (1-z^{2n-1})^5(1-z^{2n})^3
$$
The proof of this last statement is by an adaption of the proof given in [L\"o-Ro 2] when the ring $R$
is commutative (the proof is even easier in the skew-commutative case).
\bigskip
Here is Case 12:
$$
R_{12}={E(x_1,x_2,x_3,x_4,x_5)\over(x_1x_2,x_1x_3+x_2x_4+x_3x_5,x_4x_5)}
$$
and the Hilbert series is still $R_{12}(z)=1+5z+7z^2$,
but the corresponding Hilbert series for the Koszul dual $R_{12}^!$
is given by
$$
{1\over R_{12}^!(z)}=(1-2z)^2\prod_{n=1}^\infty (1-z^n)
$$
This is proved in the same way as it was proved for the corresponding $R$ in the commutative case
 [L\"o-Ro 1].
Finally here is Case 15:
$$
R_{15}={E(x_1,x_2,x_3,x_4,x_5)\over(x_1x_4+x_2x_3,x_1x_5,x_3x_4+x_2x_5)}
$$
which has $R_{15}(z) = 1+5z+7z^2$, but
$$
{1\over R_{15}^!(z)}=(1-2z)\prod_{n=1}^\infty (1-z^{2n-1})^3(1-z^{2n})^2
$$
which is proved in a similar way.
\bigskip
Now if $R$ is any of the  3 rings above,
then then Gorenstein ring  $G=R\propto I(k)$ has
Hilbert series $G(z)= 1+12z+12z^2+z^3$
and by the Gulliksen formula
$P_G(x,y)={P_R(x,y) \over 1-xy{\rm Ext}^*_R(k,R)(x,y)}$,
the B{\o}gvad formula ${\rm Ext}^*_R(k,R)(x,y)/P_R(x,y) = x^2y^2R(-{1\over xy})$ (cf. section 5 below)
and the L\"ofwall formula ${1\over P_R(x,y)} = (1+1/x)/R^!(xy)-R(-xy)/x$
we finally find that the transcendental properties of $P_G(x,y)$ are
rationally related to those of $R^!(z)$. For the validity of the B{\o}gvad formula,
cf the last section.
\bigskip
We now present the 3 other cases of $R$ with 4 quadratic relations (cases 21, 22 and 33) where the
Hilbert series is $R(z)=1+5z+6z^2$ and which lead to possibly strange Gorenstein rings and 3-manifolds $X$ with the
dimension of $H^1(X,{\bf Q})$ equal to 11. They are extremely easy to describe:
$$
R_{21} = R_{20}/(x_4x_5), \quad R_{22} = R_{20}/(x_3x_5) \quad {\rm and}\quad R_{33}=R_{15}/(x_4x_5)
$$
but the corresponding series $R^!(z)$ are unknown. But, using the Backelin et al programme BERGMAN [Ba]
we have calculated these series up to degrees 25, 14, and 33 respectively.
In the last case we found in characteristic 47 using a work-station with 48 GB of internal memory that
$$
{1\over(1-z)^2R^!_{33}(z)} = 1-3z-z^2+z^3+2z^4+3z^5+z^6+z^7-z^8-z^9-2z^{10}-z^{11}-3z^{12}-z^{13}-z^{14}
$$
$$
-z^{15}+z^{17}+z^{18}+2z^{19}+z^{20}+z^{21}+3z^{22}+z^{23}+z^{25}+z^{26}-z^{29}-z^{30}-z^{31}-z^{32}-z^{33} \ldots
$$
But to go from degree 31 to degree 32 we needed more than one week of calculations, and from degree 32 to degree 33
we needed three weeks of calculations, even with an optimal order of the variables.
But we still think that the series might be transcendental here.
\medskip
REMARK 4.1.--- All these results are in the case the characteristic of the base field is 0 (in cases 21, 22 and 33 
characteristic 47 or higher).
In Case 20 we have different $R^!_{20}(z)$ for all characteristics and the same remarks might be applicable
to the cases 21, 22 and 33. What this gives for the corresponding fundamental groups of the corresponding
3-manifolds $X$ has not been studied.
\medskip
REMARK 4.2.--- If the Yoneda ${\rm Ext}$-algebra ${\rm Ext}^*_R(k,k)$ is not finitely generated as an algebra
 then ${\rm Ext}^*_{R \propto M}(k,k)$ is not so either, since the algebra ${\rm Ext}^*_R(k,k)$ is a quotient
of  ${\rm Ext}^*_{R \propto M}(k,k)$.
One can use this for the Gorenstein ring $G_{33}=R_{33}\propto I(k)$ which has Hilbert series
$G_{33}(z) = 1+11z+11z^2+z^3$. 
Now ${\rm Ext}^*_{R_{33}}(k,k)$ needs an infinite number of generators if and only if 
the ${\rm Tor}_{3,j}^{R^!_{33}}(k,k)$ is non-zero for an infinite number of $j:s$ (cf. e.g. Theorem 3.1, (a)  in [Ro 2]).  
Indeed, by computer calculations using the ANICK command in the programme BERGMAN one obtains that
the dimensions of ${\rm Tor}_{3,j}^{R^!_{33}}(k,k)$ are 1 for $j=4$ and then $0,3,0,2,1,2$ for $j = 5,6,7,8,9,10$
and again $0,3,0,2,1,2$ for $j = 11,12,13,14,15,16$ etc. 

\mysec{5. A formula of B{\o}gvad and $m^2$-selfinjective rings.}
\bigskip
Recall that B{\o}gvad proved the following in [B{\o}]:
Let $(R,m)$ be a local commutative ring with $m^3 =0$. Assume that
$R$ has some special properties ($soc(R) = m^2$ and $R$ being a ``beast'' [B{\o}]).
 Then we have the formula for the Bass series $Bass_R(Z)$, i.e. the generating series
of ${\rm Ext}^*_R(k,R)$, the Poincar\'e-Betti series $P_R(Z)$ ( i.e. the generating
series of ${\rm Ext}^*_R(k,k)$) and $R(Z)=1+|m/m^2|Z+|m^2|Z^2$ (the Hilbert series of $R$):
$$
Bass_R(Z)/P_R(Z)= Z^2R(-{1\over Z}).
$$
In [B{\o}] this formula is broken up into two assertions, of which the first one is often valid
(the proposition ``$E(R/m^2)$'') and the other is more special.
Lescot has observed in [Les 2], [Les 4] that the proof of this in [B{\o}] boils 
down to prove that CONDITION 5.ii and CONDITION 5.iii below are valid  under some conditions:
\bigskip
CONDITION 5.i --- The natural map ${\rm Ext}^*_R(k,m)\rightarrow {\rm Ext}^*_R(k,R)$ is an epimorphism.
\medskip
CONDITION 5.ii --- The natural map ${\rm Ext}^*_R(k,m/m^2)\rightarrow {\rm Ext}^*_R(k,R/m^2)$ is an epimorphism.
\medskip
CONDITION 5.iii --- The natural map ${\rm Ext}^*_R(k,m^2)\rightarrow {\rm Ext}^*_R(k,R)$ is an epimorphism (for $*=0$ this means
that the socle of $R$ is $m^2$).
\bigskip
The CONDITION 5.i is true if $(R,m)$ is nonregular. The CONDITION 5.ii
is also often true (cf. condition ``$E(R/m^2)$'' in [B{\o}] and Proposition 1.10 in [Les 2] for 
$I=m^2$, as well as the assertion that $G \rightarrow G/soc(G)$ is a Golod
map for a Gorenstein ring $G$ [Av-Le]).
We have not yet proved all skew-commutative versions of the preceding
results, but computer computations indicate that they are true up
to ``high degrees'' and probably in all degrees.
\bigskip
We now present some results that should give the skew-versions of some of the
preceding results. We hope to return to these problems rather soon.
\medskip
Let $R$ be any ring (with unit). Recall that in order to test that a left $R$-module $M$ is injective it is
sufficient to test that for any left $R$-ideal $J$,  any $R$-module map
$\phi:  J \rightarrow M$ can be extended to a map $R \rightarrow M$. Since the last map is given as 
$r \rightarrow r\cdot m_{\phi}$ where $m_{\phi}$ is a suitable element of $M$,
this means that $\phi(j)-jm_{\phi}=0$ for any $j\in J$. In particular $R$ is self-injective
to the left if and only if for any $R$-module map $\phi: J \rightarrow R$ there is an element $r_{\phi} \in R$
such that $\phi(j)- j\cdot r_{\phi}= 0$ for all $j\in J$. This explains the condition b) in the Lemma that follows:
\bigskip
LEMMA 5.1.--- Let $(R,m)$ be a local (skew)commutative local ring where $m^3=0$ and $J$ an ideal in $R$.
 The following two conditions are equivalent:

a) The natural map
$$
{\rm Ext}^1_R(R/J,m^2) \longrightarrow {\rm Ext}^1_R(R/J,R)
$$
is surjective.

b) For any $R$-module map $\phi: J \rightarrow R$ there is an element $r_{\phi} \in R$ such that
$\phi(j)- j\cdot r_{\phi} \in m^2$ for all $j\in J$.
\bigskip
PROOF: The short exact sequence $0\rightarrow J \rightarrow R \rightarrow R/J \rightarrow 0$
and the natural map $m^2 \rightarrow R$ give rize to a commutative diagram with exact rows:
$$
 0 \rightarrow {\rm Hom}_R(R/J,m^2)\rightarrow {\rm Hom}_R(R,m^2)\rightarrow {\rm Hom}_R(J,m^2)\buildrel{\delta'}\over\rightarrow {\rm Ext}^1_R(R/J,m^2)\rightarrow  0
$$
$$
    \downarrow \quad \quad \quad \quad  \quad \quad \quad \quad \downarrow \quad \quad \quad \quad \quad \quad \quad \quad \downarrow \kappa \quad \quad \quad \quad \quad \quad \quad \quad \quad \downarrow \iota
$$
$$
0 \rightarrow {\rm Hom}_R(R/J,R)\rightarrow {\rm Hom}_R(R,R)\buildrel{i}\over\rightarrow {\rm Hom}_R(J,R)\buildrel{\delta}\over\rightarrow {\rm Ext}^1_R(R/J,R)\rightarrow  0
$$
Let us first prove that  $ a) \Rightarrow b)$. Thus assume that $\iota$ is onto. We start with
a $\phi \in {\rm Hom}_R(J,R)$ and put $\alpha = \delta(\phi)$.
We can assume that $\alpha = \iota(\xi)$. But $\xi = \delta'(\tilde\xi)$.
Thus $\delta(\kappa(\tilde\xi)-\phi)=0$, so that $\kappa(\tilde\xi)-\phi$ comes by $i$ from a map $R \rightarrow R$
of the form $ r \rightarrow r\cdot r_{\phi} $. Thus to any map $\phi : J \rightarrow R$ there is an $r_{\phi}$ in $R$
such that for all $j\in J$, $\phi(j)-j\cdot r_{\phi} \in m^2$, i.e. we have proved b). The converse follows from the ``reverse''
reasoning.
\bigskip
REMARK 5.1.--- If $L$ is any $R$-module and $P(L) \rightarrow L$ is the projective envelope of $L$ we have an exact
sequence
$$
0 \rightarrow S(L) \rightarrow P(L) \rightarrow L \rightarrow 0
$$ 
where the first syzygy $S(L)$ of $L$ is included in $m. P(L)$, we can redo the same reasoning for $P(L)/S(L)$ as we did for
$R/J$ in LEMMA 5.1. The result is that ${\rm Ext}^1_R(L,m^2) \rightarrow {\rm Ext}^1_R(L,R)$ is onto
if and only if for any map $\phi : S(L) \rightarrow R$ there is a map $j_{\phi}: P(L) \rightarrow R$ such that
$\phi(s)-j_{\phi}(s) \in m^2$ for all $s\in S(L) \subset P(L)$. Note that $P(L)$ is free so 
that $j_{\phi}$ is given by a matrix
of elements in $R$. 
\bigskip
REMARK 5.2.--- In the selfinjective case ($0$-selfinjective) it is of course sufficient to require b) for the maximal ideal $J=m$.
In the ``$m^2$-selfinjective'' case it is not clear (for us) what the right definition should be. We hope to return to this later. Therefore the definition below is
maybe too strong.
\bigskip
DEFINITION 5.1.--- We say that $R$ is $m^2$-selfinjective if the conditions of REMARK 5.2 are valid
for $L=k$ and all syzygies of $k$.
\bigskip
Using this definition we can formulate the following consequence of LEMMA 5.1 and and REMARK 5.1:
\bigskip
COROLLARY 5.1.--- The following conditions are equivalent:

$\alpha$) ${\rm Ext}^*_R(k,m^2) \rightarrow {\rm Ext}^*_R(k,R)$ is onto.
\medskip
$\beta$)  $R$ is $m^2$-selfinjective.
\bigskip
REMARK 5.3. --- In the commutative case, it seems that many rings $(R,m)$ with $m^3=0$
of the form $R=k[x_1,x_2,\ldots,x_n]/(f_1,f_2,\ldots,f_t)$ where
the $f_i$ are homogeneous quadratic forms are $m^2$-selfinjective.
For the case when $n=4$, cf. REMARK 5.4 below. 
\bigskip
OBSERVATION 5.1.--- Let $(R,m)$ be a local ring with $m^3=0$, residue field $k=R/m$ and $soc(R)=m^2$.

The following two conditions are equivalent:

a) $R$ is $m^2$-selfinjective and condition CONDITION 5.ii above is true.

b) The Bass series of $R$, i.e.
$Bass_R(Z)=\sum_{i\geq 0} |{\rm Ext}^i(k,R)| Z^i$ is related to the Poincar\'e-Betti series 
$P_R(Z) = \sum_{i\geq 0} |{\rm Ext}^i_R(k,k)| Z^i$ by the ``B{\o}gvad formula''
$Bass_R(Z) = Z^2R(-1/z)P_R(Z)$,
where $R(Z) = 1+|m/m^2|Z+|m^2|Z^2$ is the Hilbert series of $R$.

PROOF:- Consider the long exact sequence obtained when we apply ${\rm Ext}^*_R(k,-)$
to the short exact sequence $0 \rightarrow m^2 \rightarrow R \rightarrow R/m^2 \rightarrow 0$:
$$
0 \rightarrow {\rm Hom}_R(k,m^2)\rightarrow {\rm Hom}_R(k,R)\rightarrow {\rm Hom}_R(k,R/m^2)\rightarrow {\rm Ext}^1_R(k,m^2) \rightarrow {\rm Ext}^1_R(k,R)\rightarrow
$$
$$
 \rightarrow {\rm Ext}^1_R(k,R/m^2)\rightarrow {\rm Ext}^2_R(k,m^2)\rightarrow {\rm Ext}^2_R(k,R)\rightarrow {\rm Ext}^2_R(k,R/m^2)\rightarrow 
{\rm Ext}^3_R(k,m^2) \rightarrow 
$$
Since $soc(R)=m^2$ the natural monomorphism ${\rm Hom}_R(k,m^2)\rightarrow {\rm Hom}_R(k,R)$ is indeed an isomorphism.
Now according to a) all the maps ${\rm Ext}^i_R(k,m^2) \rightarrow {\rm Ext}^i_R(k,R)$ are also epimorphisms for $i \geq 1$.
It follows that we have an exact sequence:
$$
 0 \rightarrow s^{-1}{\rm Ext}^*_R(k,R/m^2)\rightarrow {\rm Ext}^*_R(k,m^2)\rightarrow {\rm Ext}^*_R(k,R)\rightarrow 0
$$
so that $Bass_R(Z)-|m^2|P_R(Z)+Z.{\rm Ext}^*_R(k,R/m^2)(Z)=0$.
But it is now easy to apply CONDITION 5.ii and this gives the result that $a)\Rightarrow b)$.
The converse is easy.
\bigskip
REMARK 5.4.--- In the commutative case most rings of embedding dimension 4, which are quotients of $k[x,y,z,u]$ by an ideal
$I$ generated by homogeneous quadratic forms and having $m^3=0$, are according to [Ro 5]
given by (the numbers of the ideals comes from [Ro 5]):
\medskip
$I_{29} = (x^2+xy,y^2+xu,z^2+xu,zu+u^2,yz)$
\medskip
$I_{54} = (x^2,xy,y^2,z^2,yu+zu,u^2)$
\medskip
$I_{55} = (x^2+xy,xz+yu,xu,y^2,z^2,zu+u^2)$
\medskip
$I_{56} = (x^2+xz+u^2,xy,xu,x^2-y^2,z^2,zu)$
\medskip
$I_{57} = (x^2+yz+u^2,xu,x^2+xy,xz+yu,zu+u^2,y^2+z^2)$
\medskip
$I_{71} = (x^2,y^2,z^2,u^2,xy,zu,yz+xu) $
\medskip
$I_{78} = (x^2,xy,y^2,z^2,zu,u^2,xz+yu,yz-xu)$
\medskip
$I_{81} = (x^2,y^2,z^2,u^2,xy,xz,yz-xu,yu,zu)$.
\medskip
The Hilbert series $R(z)$ of the different cases $R=k[x,y,z,u]/I$ are
$$
1+4z+5z^2\quad {\rm (case\,\, 29)},\quad 1+4z+4z^2\quad {\rm (cases\,\, 54, 55, 56, 57)},
$$
$$
1+4z+3z^2\quad {\rm (case\,\, 71)},\quad 1+4z+2z^2\quad {\rm (case\,\, 78)\quad and}\quad 1+4z+z^2\quad {\rm (case\,\, 81)} 
$$
But the Hilbert series $R^!(z)$ of the Koszul duals $R^!$ are respectively (all different):
$$
{(1+z)^4\over(1-z^2)^5},\quad {1\over 1-4z+4z^2},\quad {(1-z+z^2)^2\over (1-z)^3(1-3z+3z^2-3z^3)},\quad {1-z+z^2 \over (1-z)^2(1-3z+2z^2-z^3)} 
$$
$$
{1 \over (1-z)^2(1-2z-z^2)},\quad {1 \over 1-4z+3z^2},\quad  {1 \over 1-4z+2z^2}\quad{\rm and}\quad  {1 \over 1-4z+z^2}
$$
and in all cases the ``L\"ofwall formula''
$$
{1\over P_R(x,y)}= (1+{1/x})/R^!(xy)-R(-xy)/x
$$
holds true.
In all these cases {\it except} $I_{78}$ the B{\o}gvad formula also holds true.

\bigskip
REMARK 5.5.--- The B{\o}gvad formula is indeed a two-variable formula:
$$
Bass_R(x,y)/P_R(x,y)= x^2y^2R(-{1\over xy})
$$
 which shows that the non-diagonal elements occur for the two-variable
Bass series in ``the same way'' as they occur in the two-variable Poincar\'e-Betti series.
In the case of $I_{78}$ this is not true. Indeed $I_{78}$ is a Koszul ideal but the corresponding
Bass series has non-diagonal elements. More precisely we have in that case:
$$
Bass_R(x,y)/P_R(x,y) = x^2y^2R(-{1\over xy})+x^2y^2+xy^2.
$$
REMARK 5.6.--- There are of course similar results to those of REMARK 5.4 and 5.5 in the skew-commutative case.
What corresponds to the ``bad'' case $I_{78}$ in the case of four commuting variables is the case
of four skew-commuting variables: $E(x,y,z,u)\over (xy,xz-yu,yz-xu,zu)$.
\bigskip
REMARK 5.7.--- In general it is not true that $Bass_R(x,y)/P_R(x,y)$ is a rational function.
This was first noted by Lescot in his thesis [Les 4]. Take e.g. $S$ a ring with transcendental Bass series.
Form $T = S \propto I(S/m)$ and let $R$ be $T/(socle T)$ then $R$ is a so-called Teter ring [Tet], where
the maximal ideal is isomorphic to its Matlis dual (this even characterizes Teter rings [Hu]), and from
this one sees that $Bass_R(x,y)/P_R(x,y)$ is transcendental [Les 2,Corollaire 1.9 ]. 

\bigskip

APPENDIX.--- How to prove that an algebra is Koszul using non-commutative permuted Gr{\"o}bner bases, or
Macaulay2 [Gray-Sti].

\bigskip
It is well-known that a quadratic algebra (generators in degree 1, relations in degree 2) is Koszul
if it has a quadratic Gr\"obner basis [Fr\"o 2] for some ordering of the variables.
 This is stronger than being a Koszul algebra, but
most of the dual Koszul algebras of Theorem 3.2 above satisfy this stronger condition {\it for a suitable
permutation of the variables}. Since there are 9 variables, there are 9! permutations. But several years
ago J\"orgen Backelin constructed at my request a programme {\tt permutebetter.sl} written in PSL and running
under BERGMAN that in this case goes through the permutations of the variables and indicates the
length of Gr\"obner basis for each permutation of the variables.
We will only indicate how this works for the case 63 in the table of Vinberg et al.
Here the 3-form is in their notations $129\,138\,167\,246\,257\,345$ leading according to Th-Rec. 2.1 to the algebra
$$
{{\scriptstyle k<e1,e2,e3,e4,e5,e6,e7,e8,e9>}\over
{\scriptstyle ([e2,e9]+[e3,e8]+[e6,e7],[e1,e9]-[e4,e6]-[e5,e7],[e1,e8]-[e4,e5],[e2,e6]+[e3,e5],}}
$$
$$
{\scriptstyle [e2,e7]-[e3,e4],[e1,e7]-[e2,e4],[e1,e6]+[e2,e5],[e1,e3],[e1,e2])}
$$
We now construct the input file for BERGMAN (no variables mentioned -- they will be permuted) {\tt invinberg63}:

{\tts
\noindent (setq embdim 9)

\noindent (setq maxdeg 6)

\noindent (noncommify)

\noindent (off GC)

\noindent (ALGFORMINPUT)

\noindent e2*e9-e9*e2+e3*e8-e8*e3+e6*e7-e7*e6,-e1*e9+e9*e1+e4*e6-e6*e4+e5*e7-e7*e5,

\noindent -e1*e8+e8*e1+e4*e5-e5*e4,e2*e6-e6*e2+e3*e5-e5*e3,-e2*e7+e7*e2+e3*e4-e4*e3,

\noindent -e1*e7+e7*e1+e2*e4-e4*e2,e1*e6-e6*e1+e2*e5-e5*e2,e1*e3-e3*e1,e1*e2-e2*e1;
}

We also need the {\tt varsfile9} which numbers the variables e1, e2, e3, e4, e5, e6, e7, e8, e9:

\noindent {\tts ((1 . e1) (2 . e2) (3 . e3) (4 . e4) (5 . e5) (6 . e6)(7 . e7)(8 . e8)(9 . e9))}

We are now ready to go and we start BERGMAN and load in the programme needed
to do the work: {\tt permutebetter.sl} that can be found on

\noindent {\tt http://www.maths.lth.se/matematiklth/personal/ufn/bergman/permutebetter.sl}

We get a command {\tt permutedgbases} that we use as seen below
creating an outputfile {\tt outvinberg63}:

\noindent {\tts Bergman 1.001, 29-Aug-20}

\noindent {\tts 1 lisp> (dskin ``permutebetter.sl'')}

\noindent {\tts nil}

\noindent {\tts permutedgbases}

\noindent {\tts nextpermalist}

\noindent {\tts *** Function `degreeenddisplay' has been redefined}

\noindent {\tts degreeenddisplay}

\noindent {\tts nil}

\noindent {\tts nil}

\noindent {\tts 2 lisp> (permutedgbases ``varsfile9'' ``invinberg63'' ``outvinberg63'')}

The outputfile {\tt outvinberg63} contains all the information we need.
In particular we see using a PERL programme, that Torsten Ekedahl has written for me (I thank him for that), that the order of the variables: {\tt e5, e1, e2, e4, e6, e7, e3, e8, e9}
gives a non-commutative quadratic groebner basis and therefore the case 63 is Koszul.
More precisely: of the 362880 permutations there are 18618 that give a quadratic groebner basis. Indeed:

{\tts fgrep -e '\% 2' outvinberg63 | wc -l}

gives 362880, but 

{\tts fgrep -e '\% 3' outvinberg63 | wc -l}

gives 344262, and 362880-344262=18618.

However for the case 77 in the table 6 of Vinberg et al:
$129\,136\,138\,147\,234\,256$
both {\tt wc -l} give 362880 . In this case we can use Macaulay2 [Gray-Sti] on the  
skew-commutative Koszul dual and find that in this case we have a quadratic groebner basis as follows:

In this case we can use Macaulay2 on the skew-commutative Koszul dual and find in particular
that we have a quadratic groebner basis as follows (here we use {\tts ei} as a notation for the dual variables too):

{\tts
R:=QQ[e1,e2,e3,e4,e5,e6,e7,e8,e9,SkewCommutative => true]

I:=ideal(e1*e5,e1*e6,e2*e7,e2*e8,e3*e5,e3*e6,e3*e7,e3*e9,e4*e5,e4*e6,e4*e8,e4*e9,e5*e7,

e5*e8,e5*e9,e6*e7,e6*e8,e6*e9,e7*e8,e7*e9,e8*e9,e1*e7+e2*e3,e1*e8-e2*e4,

e1*e9+e3*e4,e1*e9+e5*e6,e2*e9-e3*e8,e3*e8-e4*e7)

G= gens gb I

}

The output file contains

{\tts i3 : G= gens gb I

o3 = | e8e9 e7e9 e6e9 e5e9 e4e9 e3e9 e7e8 e6e8 e5e8 e4e8 e3e8-e2e9 

e2e8 e6e7 e5e7 e4e7-e2e9 e3e7 e2e7 e5e6+e1e9 e4e6 e3e6 e1e6 e4e5 

e3e5 e1e5 e3e4+e1e9 e2e4-e1e8 e2e3+e1e7 |}

which shows that the (skew-commutative) groebner basis is quadratic and thus case 77 also gives a Koszul algebra.
\bigskip

REFERENCES
\bigskip
[An] Anick, D. J., Gulliksen, T. H.,
Rational dependence among Hilbert and Poincar\'e series. 
J. Pure Appl. Algebra 38 (1985), no. 2-3, 135--157.

[Av-Le] Avramov, L. L., Levin, G. L., 
Factoring out the socle of a Gorenstein ring. 
J. Algebra 55 (1978), no. 1, 74--83.

[Ba] Backelin, J. et al, BERGMAN, a programme for non-commutative Gr\"obner basis
calculations available at
{\tt http://servus.math.su.se/bergman/}.

[Be] Benson, Dave
An algebraic model for chains on $\Omega BG^\wedge_p$. 
Trans. Amer. Math. Soc. 361 (2009), no. 4, 2225--2242.

[B{\o}] B{\o}gvad, R.,
Gorenstein rings with transcendental Poincar\'e-series. 
Math. Scand. 53 (1983), no. 1, 5--15.

[Co-Su] Cohen, D. C., Suciu, A. I., The boundary manifold of a complex line arrangement. Groups, homotopy and configuration spaces, 105--146, Geom. Topol. Monogr., 13, Geom. Topol. Publ., Coventry, 2008.

[Eis] Eisenbud, D., Commutative algebra. With a view toward algebraic geometry.
Graduate Texts in Mathematics. 150. Berlin: Springer-Verlag, 1995. xvi+785 p.

[Elk-Sri] El Khoury, S., Srinivasan, H.,
A class of Gorenstein artin algebras of embedding dimension four. Communication in Algebra, 37 (2009), 3259-3257.

[Fo-Gr-Rei] Fossum, R. M., Griffith, P. A., Reiten, I.,
Trivial extensions of abelian categories. 
Homological algebra of trivial extensions of abelian categories with applications to ring theory. Lecture Notes in Mathematics, Vol. 456. Springer-Verlag, Berlin-New York, 1975. xi+122 pp.

[Fr\"o 1] Fr\"oberg, R.,
Determination of a class of Poincar\'e	 series. 
Math. Scand. 37 (1975), no. 1, 29--39.

[Fr\"o 2] Fr\"oberg, R.,
Koszul algebras. Advances in commutative ring theory (Fez, 1997), 337--350, 
Lecture Notes in Pure and Appl. Math., 205, Dekker, New York, 1999. 

[Gray-Sti] Grayson, D.  R., Stillman, M. E.,
Macaulay2, a software system for research in algebraic geometry.
Available at {\tt http://www.math.uiuc.edu/Macaulay2/}.

[Gu] Gulliksen, T. H.,
Massey operations and the Poincar\'e series of certain local rings. 
J. Algebra 22 (1972), 223--232.

[Gur] Gurevich, G. B.,
Foundations of the theory of algebraic invariants. 
Translated by J. R. M. Radok and A. J. M. Spencer P. Noordhoff Ltd., Groningen 1964 viii+429 pp.

[Hu] Huneke, C., Vraciu, A.,
 Rings which are almost Gorenstein.
   
{\tt arXiv:math/0403306}.

[I] Ivanov, A. F.
Homological characterization of a class of local rings. (Russian) 
Mat. Sb. (N.S.) 110(152) (1979), no. 3, 454--458, 472.

[Les 1] Lescot, J., S\'eries de Bass des modules de syzygie. [Bass series of syzygy modules] Algebra, algebraic topology and their interactions (Stockholm, 1983), 277--290, 
Lecture Notes in Math., 1183, Springer, Berlin, 1986.

[Les 2] Lescot, J.,
La s\'erie de Bass d'un produit fibr\'e	 d'anneaux locaux.  [The Bass series of a fiber product of local rings] Paul Dubreil and Marie-Paule Malliavin algebra seminar, 35th year (Paris, 1982), 218--239, 
Lecture Notes in Math., 1029, Springer, Berlin, 1983.  

[Les 3] Lescot, J.,
Asymptotic properties of Betti numbers of modules over certain rings. 
J. Pure Appl. Algebra 38 (1985), no. 2-3, 287--298.

[Les 4]  Lescot, J.,
Contribution \`a l'\'etude des s\'eries de Bass, Th\`ese, Universit\'e de Caen, 1985.

[Lev] Levin, G.L.,
Modules and Golod homomorphisms. 
J. Pure Appl. Algebra 38 (1985), no. 2-3, 299--304.

[L\"of 1] L\"ofwall, C.,
The global homological dimensions of trivial extensions of rings.
J. Algebra 39 (1976), no. 1, 287--307.

[L\"of 2] L\"ofwall, C., On the subalgebra generated by the one-dimensional elements in the Yoneda Ext-algebra. Algebra, algebraic topology and their interactions (Stockholm, 1983), 291--338, 
Lecture Notes in Math., 1183, Springer, Berlin, 1986.

[L\"of-Ro 1] L\"ofwall, C., Roos,J.-E.,
Cohomologie des alg\`ebres de Lie gradu\'ees et s\'eries de Poincar\'e-Betti non rationnelles.
 C. R. Acad. Sci. Paris S\'er. A-B 290 (1980), no. 16, A733--A736.

[L\"of-Ro 2] L\"ofwall,C., Roos,J.-E.,
A nonnilpotent $1$-$2$-presented graded Hopf algebra whose Hilbert series converges in the unit circle. 
Adv. Math. 130 (1997), no. 2, 161--200.

[Mat] Matlis, E.,
Injective modules over Noetherian rings. 
Pacific J. Math. 8 1958 511--528.

[Mil] Milnor, J.,
A unique decomposition theorem for $3$-manifolds. 
Amer. J. Math. 84 1962 1--7.

[Pos] Positselski, L. E.,
The correspondence between Hilbert series of quadratically dual algebras does not imply their having the Koszul property. (Russian) Funktsional. Anal. i Prilozhen. 29 (1995), no. 3, 83--87; translation in 
Funct. Anal. Appl. 29 (1995), no. 3, 213--217 (1996).

[Ro 1] Roos, J.-E.,
Relations between Poincar\'e-Betti series of loop spaces and of local rings. S\'eminaire d'Alg\`ebre Paul Dubreil 
31\`eme ann\'ee (Paris, 1977--1978), pp. 285--322, 
Lecture Notes in Math., 740, Springer, Berlin, 1979.

[Ro 2] Roos, J.-E.,
The homotopy Lie algebra of a complex hyperplane arrangement is not necessarily finitely presented. 
Experiment. Math. 17 (2008), no. 2, 129--143.

[Ro 3] Roos, J.-E., 
Homological properties of quotients of exterior algebras, In preparation,
Abstract available at {\it Abstracts Amer. Math. Soc.} 21 (2000), 50-51.

[Ro 4] Roos, J.-E.,
On the characterisation of Koszul algebras. Four counterexamples.
C. R. Acad. Sci. Paris S\'er. I Math. 321 (1995), no. 1, 15--20.

[Ro 5] Roos, J.-E.,
A description of the Homological Behaviour of Families of Quadratic Forms in Four variables,
Syzygies and Geometry: Boston 1995 (A.Iarrobino, \break A.Martsinkovsky and J.Weyman, editors),
Northeastern Univ. 1995, 86-95.

[Sik] Sikora, A. S.,
Cut numbers of 3-manifolds.
Trans. Amer. Math. Soc. 357 (2005), no. 5, 2007--2020 (electronic).

[Sul] Sullivan, D.,
On the intersection ring of compact three manifolds. 
Topology 14 (1975), no. 3, 275--277.

[Tet] Teter, W.,
Rings which are a factor of a Gorenstein ring by its socle. 
Invent. Math. 23 (1974), 153--162.

[Vin-El] Vinberg,B., Elashvili, A. G.,

A classification of the three-vectors of 
nine-dimensional space. (Russian)

Trudy Sem. Vektor. Tenzor. Anal. 18 (1978), 197--233.

(English translation in: Selecta Math. Soviet. 7 (1988), no. 1. 63-98.)
\par}
\end